%% file: main.tex
\documentclass[conference]{IEEEtran}
\IEEEoverridecommandlockouts
\usepackage{cite}
\usepackage{amsmath,amssymb,amsfonts}\usepackage{bm}
\usepackage{algorithm}
\usepackage{algpseudocode}
\usepackage{graphicx}
\usepackage{textcomp}
\usepackage{xcolor}
\usepackage{optidef}
\usepackage[acronym]{glossaries}
\makeglossaries

\newacronym[type=\acronymtype]{ieee}{IEEE}{Institute of Electrical and Electronics Engineers}
\newacronym[type=\acronymtype]{ml}{ML}{Machine Learning}
\newacronym[type=\acronymtype]{lmi}{LMI}{Linear Matrix Inequality}
\newacronym[type=\acronymtype]{roa}{ROA}{Region of Attraction}
\newacronym[type=\acronymtype]{ts}{TS}{Takagi-Sugeno}
\newacronym[type=\acronymtype]{lti}{LTI}{Linear Time Invariant}
\newacronym[type=\acronymtype]{ltv}{LTV}{Linear Time Varying}
\newacronym[type=\acronymtype]{lpv}{LPV}{Linear Parameter Varying}
\newacronym[type=\acronymtype]{ldi}{LDI}{Linear Differential Inclusion}
\newacronym[type=\acronymtype]{di}{DI}{Differential Inclusion}
\newacronym[type=\acronymtype]{pldi}{PLDI}{Polytopic Linear Differential Inclusion}

\newacronym[type=\acronymtype]{lp}{LP}{Linear Programming}
\newacronym[type=\acronymtype]{qp}{QP}{Quadratic Programming}
\newacronym[type=\acronymtype]{qcqp}{QCQP}{Quadratically Constrained Quadratic Programming}
\newacronym[type=\acronymtype]{socp}{SOCP}{Second Order Cone Programming}
\newacronym[type=\acronymtype]{sdp}{SDP}{Semidefinite Programming}

\newacronym[type=\acronymtype]{sos}{SOS}{Sum of Squares}
\newacronym[type=\acronymtype]{psd}{PSD}{Positive Semidefinite}
\newacronym[type=\acronymtype]{pd}{PD}{Positive Definite}

\newacronym[type=\acronymtype]{ode}{ODE}{Ordinary Differential Equations}
\newacronym[type=\acronymtype]{are}{ARE}{Algebraic Riccati Equation}
\newacronym[type=\acronymtype]{mf}{MF}{Membership Function}
\newacronym[type=\acronymtype]{mpc}{MPC}{Model Predictive Control}
\newacronym[type=\acronymtype]{mhe}{MHE}{Moving Horizon Estimation}
\newacronym[type=\acronymtype]{clf}{CLF}{Control Lyapunov Function}
\newacronym[type=\acronymtype]{cbf}{CBF}{Control Barrier Function}

\usepackage{amsthm} 

\newtheorem{theorem}{Theorem}

\newtheorem{definition}{Definition}

\newtheorem{remark}{Remark}

\newcommand{\innerProduct}[2]{\langle #1, #2 \rangle} 

\newcommand{\innerS}[1]{\langle #1 \rangle_{s}} 

\newcommand{\matrixTwoTwo}[4]{
  \begin{bmatrix}
    #1 & #2 \\
    #3 & #4
  \end{bmatrix}%
}

\newcommand{\matrixTwoOne}[2]{
  \begin{bmatrix}
    #1 \\
    #2
  \end{bmatrix}%
}

\newcommand{\rectBracket}[1]{\left[ #1 \right]} 

\newcommand{\boldCo}[1]{\mathbf{Co}\left\{#1\right\}}

\def\BibTeX{{\rm B\kern-.05em{\sc i\kern-.025em b}\kern-.08em
    T\kern-.1667em\lower.7ex\hbox{E}\kern-.125emX}}
\begin{document}

\title{Estimation of Regions of Attraction for Nonlinear Systems via Coordinate-Transformed TS Models and Piecewise Quadratic Lyapunov Functions}


\author{
\normalsize{
\begin{tabular}{ccc}
    \begin{tabular}{c}
        Artun Sel\\
        Electrical and Computer Engineering\\
        The Ohio State University\\
        Columbus, OH 43210, USA\\
        artunsel@ieee.org
    \end{tabular}
    &
    \begin{tabular}{c}
        Mehmet Koruturk\\
        Electrical and Computer Engineering\\
        Virginia Tech\\
        Blacksburg, VA 24061, USA\\
        mkoruturk@vt.edu
    \end{tabular}
    &
    \begin{tabular}{c}
        Erdi Sayar\\
        Robotics, Artificial Intelligence \\
        and Real-time Systems\\
        Technical University of Munich \\
        Munich, Germany\\ 
        erdi.sayar@tum.de
    \end{tabular}
\end{tabular}
}
}




\maketitle

\input{abstract.tex}      

\input{INTRO.tex}

\input{PRELIMINARIES.tex}


\input{ROA_TS.tex}      

\input{ROA_CoC.tex}      

\input{CONCLUSIONS.tex}      




\bibliographystyle{IEEEtran}
\bibliography{Ref.bib}

\end{document}

%% file: abstract.tex
\begin{abstract}
This paper presents a novel approach for computing enlarged \glspl{roa} for nonlinear dynamical systems through the integration of multiple coordinate transformations and piecewise quadratic Lyapunov functions within the \gls{ts} modeling framework. While existing methods typically follow a single-path approach of original system $\rightarrow$ \gls{ts} model $\rightarrow$ \gls{roa} computation, the proposed methodology systematically applies a sequence of coordinate transformations to generate multiple system representations, each yielding distinct \gls{roa} estimations. Specifically, the approach transforms the original nonlinear system using transformation matrices $T_1, T_2, \ldots, T_N$ to obtain $N$ different coordinate representations, constructs corresponding \gls{ts} models for each transformed system, and computes individual \glspl{roa} using piecewise quadratic Lyapunov functions. The final \gls{roa} estimate is obtained as the union of all computed regions, leveraging the flexibility inherent in piecewise quadratic Lyapunov functions compared to traditional quadratic approaches. The enhanced methodology demonstrates significant improvements in \gls{roa} size estimation compared to conventional single-transformation techniques, as evidenced through comparative analysis with existing \gls{ts}-based stability methods. 
\end{abstract}


%% file: INTRO.tex
\section{Introduction}  \label{sec:Introduction}

The analysis and control of nonlinear dynamical systems constitute a cornerstone challenge in contemporary control theory \cite{Kurkcu18:Disturbance,Kurkcu18:DisturbanceUncertainty,Kurkcu19:Robust}. 
This challenge has gained heightened significance with the proliferation of autonomous aerial vehicles and their widespread deployment in remote sensing applications \cite{Hind22:Assessment}, the emergence of critical cybersecurity vulnerabilities in communication networks \cite{Etcibasi24:Coverage,Mohammadi25:Detection,Mohammadi25:GPS,Mohammadi25:DetectionICAIC,Ahmari25:AData,Ahmari25:Evaluating}, and the increasing complexity of modern energy systems \cite{Agrawal25:Performance,Himmelstoss25:Floating,Koch25:Qualitative,Albasheri25:Energy,Ekawita25:Investigating,Pham25:Potential,kopur25:ARSFQ}. Consequently, safety-critical considerations have become paramount in system design and operation. Furthermore, stability analysis has emerged as a crucial concern in advanced computational models \cite{Sel25:LLMs,Sel24:Skin,Sel24:Algorithm} and sophisticated optimization algorithms \cite{Gu25:Safe,Sel23:Learning,Khattar23:CMDPWithinOnline,Coskun22:Magnetic}, particularly within the context of state estimation problems critical to aerospace applications \cite{Pourtakdoust23:Advanced,Nasihati22:Satellite,Nasihati21:OnLine,Nasihati19:Autonomous}.

The analysis of stability properties and \gls{roa} estimation for nonlinear dynamical systems represents one of the fundamental challenges in control theory and applications. For a given nonlinear system of the form $\dot{x} = f(x)$, determining the largest set of initial conditions from which trajectories converge to an equilibrium point is crucial for understanding system behavior and designing robust controllers. However, exact computation of \glspl{roa} for general nonlinear systems remains analytically intractable, necessitating the development of systematic approximation methods \cite{Tedrake10:LQR}.

Among the various approaches for nonlinear system analysis, \gls{ts} fuzzy models have emerged as a powerful framework due to their ability to exactly represent nonlinear dynamics within compact domains while maintaining a linear-like structure that facilitates stability analysis \cite{LopezEstrada17:LMI_based_fault}. The \gls{ts} modeling approach employs sector nonlinearity techniques to construct a convex combination of linear subsystems, enabling the application of \gls{lmi} based methods for stability verification and \gls{roa} estimation.

Traditional approaches for \gls{roa} computation using \gls{ts} models typically rely on quadratic Lyapunov functions, which, while computationally tractable, often yield conservative estimates due to their inherent structural limitations. To address this conservatism, several researchers have proposed the use of non-quadratic Lyapunov functions, including fuzzy Lyapunov functions \cite{Bernal2022:BOOK} and piecewise quadratic Lyapunov functions \cite{DellaRossa20:MaxMinLyapunov}. The latter approach has demonstrated particular promise in reducing conservatism while maintaining computational feasibility through piecewise linear partitioning of the state space.

Despite these advances, existing methodologies typically follow a single-pathway approach: original nonlinear system $\rightarrow$ \gls{ts} model construction $\rightarrow$ \gls{roa} computation. This conventional paradigm inherently limits the achievable \gls{roa} size, as the choice of coordinates and the resulting \gls{ts} model structure significantly influence the final stability domain estimation. The fundamental insight underlying this limitation is that different coordinate representations of the same nonlinear system can yield substantially different \gls{ts} models, each potentially capturing distinct aspects of the system's stability properties.

Coordinate transformations have been extensively studied in nonlinear control theory for various purposes, including feedback linearization \cite{Khalil15:Nonlinear_Control}, normal form transformations, and system simplification. However, their systematic exploitation for \gls{roa} enlargement in the context of \gls{ts} modeling has received limited attention. The potential of coordinate transformations lies in their ability to expose different geometric structures of the nonlinear system, which can be more favorably captured by specific \gls{ts} model constructions \cite{Sotiropoulos18:Causality,Igarashi20:ARobust}.

This paper addresses the aforementioned limitations by proposing a novel multi-transformation approach that systematically leverages multiple coordinate representations to enhance \gls{roa} estimation. The key innovation lies in generating a family of transformed systems using carefully selected transformation matrices $T_1, T_2, \ldots, T_N$, constructing corresponding \gls{ts} models for each representation, and computing individual \glspl{roa} using advanced piecewise quadratic Lyapunov functions. The final \gls{roa} estimate is obtained as the union of all computed regions, effectively exploiting the complementary information provided by different coordinate perspectives.

The main contributions of this work are threefold: (i) development of a systematic framework for multiple coordinate transformation-based \gls{roa} enhancement, (ii) integration of piecewise quadratic Lyapunov functions within the multi-transformation paradigm to further reduce conservatism, and (iii) demonstration of significant \gls{roa} enlargement compared to conventional single-transformation approaches through comprehensive numerical examples.

The remainder of this paper is organized as follows. Section~\ref{sec:Preliminaries} provides the necessary background on \gls{ts} modeling and piecewise quadratic Lyapunov functions. Section~\ref{sec:ROA_TS} presents the proposed \gls{roa} estimation method utilizing \gls{ts} modeling framework with piecewise quadratic Lyapunov functions. Section~\ref{sec:ROA_CoC} introduces the novel multi-transformation methodology and demonstrates its superiority over the conventional single-transformation approach (presented in Section~\ref{sec:ROA_TS}) through comprehensive comparative analysis. Section~\ref{sec:CONCLUSIONS} concludes the paper with a summary of key findings and discussions on future research directions.

%% file: PRELIMINARIES.tex
\section{Preliminaries}  \label{sec:Preliminaries}

\subsection{\textbf{Definitions for Mathematical Optimization Related Concepts}}
\subsubsection{\textbf{Linear Matrix Inequality}}
An \gls{lmi} is a constraint of the form:
\begin{equation}
F(x) = F_0 + \sum_{i=1}^m x_i F_i \succeq 0
\end{equation}
where $x = (x_1, x_2, \ldots, x_m)^T$ represents the decision variable vector, $F_0, F_1, \ldots, F_m \in \mathbb{R}^{n \times n}$ denote given symmetric matrices, and the inequality symbol $\succeq 0$ indicates positive semidefiniteness. A \gls{lmi} inherently defines a convex constraint on the variable vector $x$. When multiple \glspl{lmi} are present, they can be consolidated into a single equivalent \gls{lmi} through the utilization of block-diagonal matrix structures. \glspl{lmi} constitute fundamental components in \gls{sdp} frameworks, where they function as convex constraints in optimization problems that admit efficient numerical solutions via interior-point algorithms \cite{boyd2004convex}.


\subsubsection{\textbf{Semi Definite Programming}}
\gls{sdp} is a powerful optimization framework that extends linear programming to the cone of \gls{psd} matrices. Formally, an \gls{sdp} can be expressed as:
\begin{equation}
\begin{aligned}
\text{min.} \quad & \innerProduct{C}{X} \\
\text{s.t.} \quad & \innerProduct{A_i}{X} \leq b_i, \quad i = 1, 2, \ldots, m \\
& X \succeq 0
\end{aligned}
\end{equation}
\noindent where $X, C, A_1, \ldots, A_m$ are symmetric $n \times n$ matrices, $b_1, \ldots, b_m$ are scalars, and $X \succeq 0$ denotes that $X$ is \gls{psd} \cite{Tedrake10:LQR}.




\subsection{\textbf{Definitions for Applications of the Methods}}

\subsubsection{\textbf{Dynamical System}}
A dynamical system constitutes a mathematical framework for modeling the temporal evolution of a physical system's state variables according to deterministic evolution laws. Formally, a dynamical system is characterized by a state space $\mathcal{X}$ and a mapping $f: \mathcal{X} \times \mathbb{R} \to \mathcal{X}$ that governs the evolution of the system state $x(t) \in \mathcal{X}$ as a function of time. In the continuous-time domain, this temporal evolution is mathematically represented by the differential equation:
\begin{equation}
\dot{x}(t) = f(x(t), t)
\end{equation}
\noindent where $\dot{x}(t)$ denotes the time derivative of the state vector. 

Dynamical systems are systematically categorized according to several fundamental characteristics: linearity versus nonlinearity, time-invariance versus time-variance, and autonomy versus non-autonomy. The theoretical analysis of such systems encompasses the investigation of stability properties, characterization of equilibrium manifolds, identification of limit cycles, and determination of other qualitative dynamical behaviors that govern system performance \cite{Duan13:LMIs}.

\begin{remark}
    In this study, we primarily focus on polynomial-type nonlinear systems that are time-invariant. For analysis problems, we specifically examine autonomous systems that possess these properties.
\end{remark}

\subsubsection{\textbf{Stability Definitions}} 
Consider a nonlinear dynamical system described by:
\begin{equation}
\label{eq:autonomous_NL_DYN_SYS}
\dot{x} = f(x),\quad x(0) = x_0
\end{equation}
where $x \in \mathbb{R}^n$ is the state vector, $f: \mathbb{R}^n \rightarrow \mathbb{R}^n$ is a locally Lipschitz function, and $f(0) = 0$ (i.e., the origin is an equilibrium point). We define various notions of stability as follows:

\begin{definition}[\textbf{Stability in the Sense of Lyapunov}]
The equilibrium point $x = 0$ is stable in the sense of Lyapunov if, for any $\epsilon > 0$, there exists a $\delta > 0$ such that:
\begin{equation}
\|x(0)\| < \delta \implies \|x(t)\| < \epsilon, \quad \forall t \geq 0
\end{equation}
\end{definition}

\begin{definition}[\textbf{Asymptotic Stability}]
The equilibrium point $x = 0$ is asymptotically stable if it is stable in the sense of Lyapunov and there exists a $\delta > 0$ such that:
\begin{equation}
\|x(0)\| < \delta \implies \lim_{t \rightarrow \infty} \|x(t)\| = 0
\end{equation}
\end{definition}


\begin{definition}[\textbf{Global Stability}]
A stability property (Lyapunov, asymptotic, or exponential) is said to be global if it holds for any initial state $x(0) \in \mathbb{R}^n$ (i.e., $\delta \rightarrow \infty$).
\end{definition}

\begin{definition}[\textbf{Local Stability}]
A stability property is said to be local if it holds only for initial states within some bounded neighborhood of the equilibrium (i.e., $\|x(0)\| < \delta$ for some finite $\delta > 0$).
\end{definition}

\begin{definition}[\textbf{\gls{roa}}]
For a given locally asymptotically stable dynamical system, the corresponding region around the equilibrium point is called \gls{roa}.
\end{definition}
\subsubsection{\textbf{Lyapunov Function and Stability Conditions for Nonlinear Dynamical Systems}}
This section first introduces the Lyapunov function, then for certain dynamical systems, a set of stability conditions (that are numerically tractable) are derived by using Lyapunov functions.   
\begin{theorem}
For \eqref{eq:autonomous_NL_DYN_SYS}, where the origin is the equilibrium point,
If $\exists V(x): \mathcal{D} \rightarrow \mathbb{R}$ \textit{continuously differentiable} such that
\begin{align*}
    V(0) &= 0, \\
    \dot{V}(0) &= 0, \\
    V(x) &> 0, \quad \forall x \in \mathcal{D} \setminus \{0\}, \\
    \dot{V}(x) &< 0, \quad \forall x \in \mathcal{D} \setminus \{0\},
\end{align*}
then $\mathbf{x} = 0$ is \textit{asymptotically stable}; if, in addition, $\mathcal{D} = \mathbb{R}^n$ and $V(x)$ is \textit{radially unbounded}, then $\mathbf{x} = 0$ is \textit{globally asymptotically stable}.
If $V(x)$ satisfies the constraints in $\Omega \subseteq \mathcal{D}$, then $\mathbf{x} = 0$ is \textit{locally asymptotically stable} and $\Omega$ is called \textit{\gls{roa}}.
\end{theorem}

\subsubsection{\textbf{Analysis Problem in \gls{lti} systems}} 
Consider a continuous-time \gls{lti} system:
\begin{equation}
\dot{x}(t) = Ax(t)
\end{equation}
where $x(t) \in \mathbb{R}^n$ is the state vector and $A \in \mathbb{R}^{n \times n}$ is the system matrix. By using a quadratic Lyapunov function parameterized by $P$, the stability condition can be written as
\begin{equation}
\begin{aligned}
\text{find} \quad & P \\
\text{s.t.} \quad & P  \in \mathbb{S}_{++} \\
& -\innerS{P A}  \in \mathbb{S}_{++}
\end{aligned}    
\end{equation}
which is a convex feasibility problem.
\subsubsection{\textbf{Quadratic Stability for parameter varying systems}}
For a given continuous-time linear system:
\begin{equation}
\dot{x}(t) = A(p)x(t)
\end{equation}
where $x(t) \in \mathbb{R}^n$ is the state vector and $A(p) \in \mathbb{R}^{n \times n},\forall p \in \mathcal{P}$ is the system matrix where $p$ is a parameter of the dynamical system, and that representation is quite general since uncertain, parameter varying, switching and differential inclusion type systems can also be represented in that way. 

Using a quadratic Lyapunov function, the quadratic stability condition amounts to the following feasibility problem:
\begin{equation}
\begin{aligned}
\text{find} \quad & P \\
\text{s.t.} \quad & P  \in \mathbb{S}_{++} \\
& -\innerS{P A(p)}  \in \mathbb{S}_{++},\forall p \in \mathcal{P}
\end{aligned}    
\end{equation}
If the set $\mathcal{P}$ is finite, the given problem is tractable and if it is not, there are methods to derive tractable sufficient conditions for stability.

\subsubsection{\textbf{Stability conditions for \gls{ldi} using quadratic and piecewise-quadratic Lyapunov functions}} \label{sec:PWC_Lyap_sec}
Consider a continuous-time dynamical system:
\begin{equation} \label{eq:def_general_DI}
\dot{x}(t) \in \boldCo{f_i(x)}, \forall i \in \rectBracket{M}
\end{equation}
where $x(t) \in \mathbb{R}^n$ is the state vector and $\boldCo{\cdot}$ denotes the convex hull of a given set. The system given in \eqref{eq:def_general_DI} is called \gls{di}, a special case where $f_i(x)=A_i x$ is called \gls{ldi}. 

Using a quadratic Lyapunov function, the quadratic stability condition amounts to the following feasibility problem:
\begin{equation}
\begin{aligned}
\text{find} \quad & P \\
\text{s.t.} \quad & P  \in \mathbb{S}_{++} \\
& -\innerS{P A_i}  \in \mathbb{S}_{++},\forall i \in \rectBracket{M}
\end{aligned}    
\end{equation}
It is important to note that the given condition is conservative, i.e., there are LDI systems that are not quadratically stable by asymptotically stable. The following condition is given to address this issue.

A piecewise quadratic Lyapunov function is given by $V(x)=\max_{j\in N}\{ x^{\top}P_jx \}$ and offers more flexibility when it comes to stability analysis. Since the stability condition for a general case is involved, the stability condition for a given LDI using a piecewise quadratic Lyapunov function with $N=2$ is given by:
\begin{equation} \label{eq:OPT_p1p2}
\begin{aligned}
\text{find} \quad & (P_1,P_2) \\
\text{s.t.} \quad & P_1,P_2 \in {{\mathbb{S}}_{++}} \\
&\exists \lambda_i :  \lambda_i (P_2-P_1)-\innerS{P_1 A_i} \in \mathbb{S}_{++},\forall i \in \rectBracket{M} \\
&\exists \lambda_i :  \lambda_i (P_1-P_2)-\innerS{P_2 A_i} \in \mathbb{S}_{++},\forall i \in \rectBracket{M}
\end{aligned}    
\end{equation}


\subsubsection{\textbf{Convex Optimization Based Nonlinear Dynamical System Analysis}} \label{sec:TS_sec}
\begin{definition}
A \textit{convex model for autonomous systems} is defined as a collection of first-order \glspl{ode} in which the right-hand side vector field can be formulated as a convex combination of constituent vector fields (sub-models), specifically
\begin{equation}
\label{eq:convex_model_autonomous}
    \dot{\mathbf{x}}(t) = \sum_{i=1}^{r} h_i(\mathbf{z}) f_i(\mathbf{x})
\end{equation}
where $r \in \mathbb{N}$ denotes the number of linear models comprising the convex combination, $\mathbf{z} \in \mathbb{R}^p$ represents a premise variable that may be a function of the system state $\mathbf{x}$, time $t$, external disturbances, or exogenous parameters. The nonlinear weighting functions $h_i(\mathbf{z})$ for $i \in {1, 2, \ldots, r}$ satisfy the \textit{convex sum property} within a compact set $\mathcal{C} \subseteq \mathbb{R}^p$, specifically:
\begin{equation}
\label{eq:convex_sum_property}
    \sum_{i=1}^{r} h_i(\mathbf{z}) = 1, \quad 0 \leq h_i(\mathbf{z}) \leq 1, \quad \forall \mathbf{z} \in \mathcal{C}.
\end{equation}
and the special case is given by
\begin{equation}
\label{eq:convex_sum_Ai_autonomous}
    \dot{\mathbf{x}}(t) = \sum_{i=1}^{r} h_i(\mathbf{z}) A_i \mathbf{x}(t)
\end{equation}
\end{definition}
%

%
\begin{definition}
The vector $\mathbf{z}$ in a convex model~\eqref{eq:convex_sum_Ai_autonomous} is known as the \textit{premise} or \textit{scheduling vector}; it is assumed to be bounded and continuously differentiable in a compact set $\mathcal{C} \subseteq \mathbb{R}^p$ of the scheduling/premise space.
\end{definition}

\begin{definition}
The nonlinear functions $h_i(\mathbf{z}),\ i \in \{1, 2, \ldots, r\}$ in a convex model~\eqref{eq:convex_sum_Ai_autonomous} are known as \textit{\glspl{mf}}; they hold the convex sum property in a compact set $\mathcal{C} \subseteq \mathbb{R}^p$ of the scheduling space.
\end{definition}

\begin{definition}
If $f_i(\mathbf{x})$ are linear, i.e., $f_i(\mathbf{x}) = A_i \mathbf{x}$, \eqref{eq:convex_sum_Ai_autonomous} is called a \textit{\gls{ts}} model~\cite{Bernal2022:BOOK}; if, in addition, $\mathbf{z}$ does not include the state $\mathbf{x}$ it is called a \textit{linear polytopic model}~\cite{Takagi1985:FuzzyIdentification}; if the \gls{ts} model comes from fuzzy modeling techniques it is also referred to as \textit{\gls{ts} fuzzy model}.
\end{definition}



\begin{remark}
When the convex sum property is satisfied globally, that is, when $\mathbf{z} \in \mathcal{C}$ for all $\mathbf{x}$ in the domain of interest, the stability characteristics established for the corresponding \gls{ldi} system—derived from the \gls{ts} representation of the original nonlinear dynamical system—directly guarantee the stability of the underlying nonlinear system.
\end{remark}

%% file: ROA_TS.tex
\section{\gls{roa} Computation Using \gls{ts} Method}  \label{sec:ROA_TS}
In this section, the \gls{roa} estimation using \gls{ts} method is described by means of a numerical example.

Consider the dynamical system described by the following equations
\begin{equation}
\begin{aligned}
    \dot{x}_1 &= \ -(x_1)^2 -2(x_2)-2(x_1) \\
    \dot{x}_2 &= \ (x_2)^3-x_2
\end{aligned}
\end{equation}
where an estimate of the \gls{roa} needs to be computed by using \gls{ts}-convex modeling.

\begin{figure}[h!]
\vspace{-0.5cm}
	\centering
	\includegraphics[width=\linewidth]{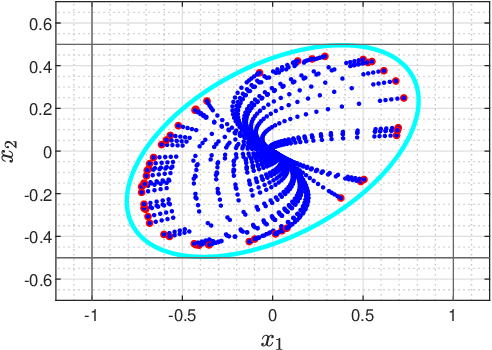}
	\caption{\gls{roa} for the system computed with the \gls{ts} method}
	\label{fig:RoA_problem_1}
\end{figure}

Considering the region given by
\begin{equation}
    \mathcal{D} = \left\{ x \in \mathbb{R}^{2} : x_1 \in [-1,+1] , x_2 \in [-0.5,+0.5] \right\}
\end{equation}
the dynamics can be rewritten as
\begin{equation}
    \matrixTwoOne{\dot{x}_1}{\dot{x}_2} = \matrixTwoTwo{-x_1-2}{-2}{0}{(x_2)^2-1)}
\matrixTwoOne{x_1}{x_2}
\end{equation}
and this can be written as
\begin{equation}
\begin{aligned}
\dot{x}(t) &= \sum_{i_1=0}^{1} \sum_{i_2=0}^{1}  
w^{1}_{i_1}(z) w^{2}_{i_2}(z)
\left(
\matrixTwoTwo{-z_{1}^{i_1}-2}{-2}{0}{z_{2}^{i_2}-1)}
\matrixTwoOne{x_1}{x_2}
\right) \\
&= \sum_{\mathbf{i} \in \mathbb{B}^2} \mathbf{w}_{\mathbf{i}}(z) 
\left( A_{\mathbf{i}} x(t)\right)
= A_{\mathbf{w}} x(t),
\end{aligned}
\end{equation}
where the terms are defined as
\begin{equation}
    z_{1}=x_1, \quad z_{2}=(x_2)^2.
\end{equation}
and are defined on the set given by
\begin{equation}
    \mathcal{C}=
    \left\{
    z \in \mathbb{R}^2 : \exists x\in \mathcal{D} \, s.t. \, z_1=x_1,z_2=(x_2)^2   
    \right\}.
\end{equation}
Explicitly, the terms are written as
\begin{equation}
        z_{1} = \sum_{i_1=0}^{1} w^{1}_{i_1}(x) z_{1}^{i_1},\quad z_{2} = \sum_{i_2=0}^{1} w^{2}_{i_2}(x) z_{2}^{i_2}.
\end{equation}
%
where the boundary-terms are computed by
\begin{equation}
    \begin{aligned}
    &\min_{x_1(t),\, x_2(t)} z_1(t) = z_{1}^{i_1=0}, \quad && \max_{x_1(t),\, x_2(t)} z_1(t) = z_{1}^{i_1=1}, \\
    &\min_{x_1(t),\, x_2(t)} z_2(t) = z_{2}^{i_2=0}, \quad && \max_{x_1(t),\, x_2(t)} z_2(t) = z_{2}^{i_2=1}.
    \end{aligned}    
\end{equation}
and the corresponding \glspl{mf} are given by
\begin{equation}
    \begin{aligned}
        w_{i_1=1}^{1}(x)=\frac{z_{1}^{1} - z_1}{z_{1}^{1}-z_{1}^{0}} &, \quad w_{i_1=0}^{1}(x)=1-w_{i_1=1}^{1}(x). \\
        w_{i_2=1}^{2}(x)=\frac{z_{2}^{1} - z_2}{z_{2}^{1}-z_{2}^{0}} &, \quad w_{i_2=0}^{2}(x)=1-w_{i_2=1}^{2}(x).    
    \end{aligned}
\end{equation}
Using the \gls{ts} method and piecewise Lyapunov function, described in Section~\ref{sec:TS_sec} and Section~\ref{sec:PWC_Lyap_sec}, respectively, the problem can be stated by the feasibility problem given in \eqref{eq:OPT_p1p2} where the state matrices are given by
\begin{equation*}
\begin{aligned}
A_1=\matrixTwoTwo{-1}{-2}{0}{-1} &, A_2=\matrixTwoTwo{-1}{-2}{0}{-0.75}, \\
A_3=\matrixTwoTwo{-3}{-2}{0}{-1} &, A_4=\matrixTwoTwo{-3}{-2}{0}{-0.75}.
\end{aligned}
\end{equation*}
and the $(P_1,P_2)$ terms are given by
\begin{equation*}
    P_1=\matrixTwoTwo{0.1071}{-0.0829}{*}{0.2836},
    P_2=\matrixTwoTwo{0.1045}{-0.0852}{*}{0.2605}.
\end{equation*}
Since the \gls{ts}-model is valid in $\mathcal{D}$, using the parameters of the Lyapunov function $(P_1,P_2)$, the \gls{roa} is computed by the following optimization problem:
\begin{equation}
    \begin{aligned}
    \label{eq:k_for_RoA_ALGO}
    \text{max.} \quad & k \\
    \text{s.t.} \quad & \Omega=\left\{ x \in \mathcal{D}: V(x) \leq k \right\} \\
                     & \Omega \subseteq \mathcal{D}
    \end{aligned}
\end{equation}
and $k$ is computed to be $0.054$
using the method described in \cite{Robles17:Subspace_Based_TS}
and the corresponding \gls{roa} is illustrated in Fig-\ref{fig:RoA_problem_1}.

%% file: ROA_CoC.tex
\section{\gls{roa} Computation Using Combined State Transformation and \gls{ts} Method}  \label{sec:ROA_CoC}
Regarding the stability analysis methodology, the \gls{ts} approach provides a systematic and computationally efficient algorithm for constructing Lyapunov functions that ensure local asymptotic stability guarantees. However, the efficacy of this method is fundamentally dependent on the selection of \glspl{mf}, and as previously established, this selection process presents significant challenges due to the non-uniqueness of MF sets for any given dynamical system.

Rather than pursuing the optimization of MF selection to maximize the \gls{roa}, this work proposes an alternative strategy. We employ a computationally tractable methodology for determining an appropriate set of \glspl{mf} for a given dynamical system, which yields an initial \gls{roa} estimate. Subsequently, through the implementation of a coordinate transformation $\bar{x} = Tx$ that results in the transformed system dynamics $\dot{\bar{x}} = \bar{f}(\bar{x})$, and applying the identical \gls{ts} methodology to this transformed system, an additional \gls{roa} can be computed.

\begin{figure}[h!]
\vspace{-1.5em}
	\centering
	\includegraphics[width=\linewidth]{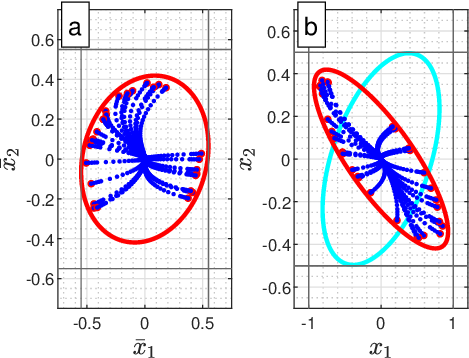}
	\caption{(a) \gls{roa} computed for the system whose state variable is $\bar{x}\in\bar{\mathcal{D}}$
    (b) \gls{roa} computed for the system whose state variable is $x\in\mathcal{D}$ by using the \gls{ts}-method illustrated by red-shaded region and the boundary of the \gls{roa} that is computed in $\bar{\mathcal{D}}$ domain and mapped to $\mathcal{D}$}
	\label{fig:RoA_problem_3}
\end{figure}

By choosing a simple transformation matrix $T=\matrixTwoTwo{1}{2}{0}{1}$,
which results in the following system:
\begin{equation}
\begin{aligned}
    \dot{\bar{x}}_1 &= \ -(\bar{x}_1^2) +4(\bar{x}_1)(\bar{x}_2)-2(\bar{x}_1)+2(\bar{x}_2^3)-4(\bar{x}_2^2)  \\
    \dot{\bar{x}}_2 &= \ (\bar{x}_2^3)-\bar{x}_2 
\end{aligned}
\end{equation}
for the following set of \glspl{mf}: $z_1 = \bar{x}_1,z_2 = \bar{x}_2,z_3 = \bar{x}_2^2$
and for the set:
\begin{equation}
\bar{\mathcal{D}}=\left\{  x\in\mathbb{R}^{2}:\bar{x}_1,\bar{x}_2\in[-0.55,+0.55] \right\},
\end{equation}
the problem can be stated by the feasibility problem given in \eqref{eq:OPT_p1p2} where the state matrices are given by
\begin{equation*}
\begin{aligned}
A_1=\matrixTwoTwo{-1.45}{-0.3328}{0}{-1.1664}&, A_2=\matrixTwoTwo{-1.45}{0.3328}{0}{-0.8336}, \\
A_3=\matrixTwoTwo{-0.24}{-1.5428}{0}{-1.1664}&, A_4=\matrixTwoTwo{-0.24}{0.8773}{0}{-0.8336}, \\
A_5=\matrixTwoTwo{-2.55}{-0.3328}{0}{-1.1664}&, A_6=\matrixTwoTwo{-2.55}{0.3328}{0}{-0.8336}, \\
A_7=\matrixTwoTwo{-1.34}{-1.5428}{0}{-1.1664}&, A_8=\matrixTwoTwo{-1.34}{-0.8773}{0}{-0.8336}.
\end{aligned}
\end{equation*}
and the $(P_1,P_2)$ terms are given by
\begin{equation*}
    P_1=\matrixTwoTwo{5.0473}{-1.1747}{*}{8.4518},
    P_2=\matrixTwoTwo{5.0896}{-1.0599}{*}{8.7648}.
\end{equation*}
the optimization problem for the given $\bar{\mathcal{D}}$ set is computed using \eqref{eq:k_for_RoA_ALGO}, resulting in $k=1.54$ and the \gls{roa} for the $\bar{x}$-domain and the region that is computed by mapping that region back to the original $x$-domain is illustrated in Fig-\ref{fig:RoA_problem_3}.

The ellipsoidal boundary depicted by the red ellipse, derived through the proposed change of variables methodology, exhibits regions that extend beyond the original \gls{roa} computed in Section \ref{sec:ROA_TS}. This observation indicates that the iterative application of the presented approach, achieved through the systematic automation of monomial function selection for a given nonlinear dynamical system, yields an enlarged \gls{roa} compared to conventional methods.

%% file: CONCLUSIONS.tex
\section{Conclusions}  \label{sec:CONCLUSIONS}
In this study, 
the \gls{ts} modeling framework was employed to study the \gls{roa} estimation problem. The \gls{ts} approach enabled the leveraging of convex optimization techniques to handle nonlinear dynamics effectively. Furthermore, an extension to the standard methods was proposed, aimed at enhancing their applicability and reducing conservatism in practical scenarios. These contributions collectively demonstrate the potential of combining \gls{ts} modeling with \gls{sdp}-based tools to address challenging problems in nonlinear dynamical systems.